\documentclass[letterpaper, 10 pt, conference]{ieeeconf}  

\IEEEoverridecommandlockouts                              
\usepackage{amsfonts, amsmath, bm}
\usepackage{algorithmic, graphicx}
\usepackage{epstopdf}
\usepackage{float}
\usepackage{subcaption}
\usepackage{color}
\usepackage{xcolor}

\title{\LARGE \bf
Robust Quantum Gate Preparation in Open Environments
} 
\author{Luke S. Baker$^1$, Syed A. Shah$^1$, Anatoly Zlotnik$^3$, Andrei Piryatinski$^4$
\thanks{The authors are grateful to Andrew Harter and Andre Luiz P. de Lima for valuable discussions. This project was supported by the LDRD program and the Center for Nonlinear Studies at Los Alamos National Laboratory.  Research conducted at Los Alamos National Laboratory is done under the auspices of the National Nuclear Security Administration of the U.S. Department of Energy under Contract No. 89233218CNA000001. Report number: LA-UR-24-30519.}
\thanks{$^1$\{lsbaker,shah\}@lanl.gov, \,\,  Center for Nonlinear Studies, Theoretical Division, Los Alamos National Laboratory, Los Alamos, NM, USA 87545}
\thanks{$^3$azlotnik@lanl.gov, \,\, Applied Mathematics \& Plasma Physics, Los Alamos National Laboratory, Los Alamos, NM, USA 87545}
\thanks{$^4$apiryat@lanl.gov, \,\, Physics of Condensed Matter \& Complex Systems, Los Alamos National Laboratory, Los Alamos, NM, USA 87545}}

\newtheorem{thm}{Theorem}

\newtheorem{prop}[thm]{Proposition}

\begin{document}

\maketitle
\thispagestyle{empty}
\pagestyle{empty}

\begin{abstract}
We develop an optimal control algorithm for robust quantum gate preparation in open environments with the state of the quantum system represented using the Lindblad master equation.  The algorithm is based on adaptive linearization and iterative quadratic programming to progressively shape the control signal into an optimal form.  Robustness is achieved with exponential rates of convergence by introducing uncertain parameters into the master equation and expanding the parameterized state over the basis of Legendre polynomials.  We prove that the proposed control algorithm reduces to GRadient Ascent Pulse Engineering (GRAPE) \cite{khaneja2005optimal} when the robustness portion of the algorithm is bypassed and signal restrictions are relaxed. The control algorithm is applied to prepare Controlled NOT and SWAP gates with high precision. Using only second order Legendre polynomials, the examples showcase unprecedented robustness to 100\% parameter uncertainty in the interaction strength between the qubits, while simultaneously compensating for 20\% uncertainty in signal intensity.  The results could enable new capabilities for robust implementation of quantum gates and circuits subject to harsh environments and hardware limitations.
\end{abstract}

\section{Introduction} \label{sec:intro}

High-performance quantum gates are expected to play a critical role in the realization of quantum computers \cite{williams2011quantum}.  The ability to achieve quantum advantage in applications of practical interest is, however, limited by the complexity of gate preparation \cite{pareek2024demystifying}.  This challenge has compelled the co-design approaches that use classical methods to enable more complex quantum computations.  
Optimal control has become one of the preferred methods over bang-bang and Lyapunov techniques for the control of quantum systems because of its flexibility and inherent optimality \cite{koch2022quantum,morton2006bang,wen2016preparation}.  A gradient-based algorithm known as GRadient Ascent Pulse Engineering (GRAPE) \cite{khaneja2005optimal} is an established method that was originally designed for optimal control in nuclear magnetic resonance (NMR) spectroscopy. This algorithm synthesizes piece-wise constant control signals in which the amplitudes are treated as optimization variables.  

Since the time of its development, numerous extensions have emerged to  advance the functionality of GRAPE as well as to accommodate broader classes of quantum systems \cite{de2011second,motzoi2011optimal, bhole2018practical, chen2023accelerating}.  The Krotov method \cite{tannor1992control, eitan2011optimal} offers an alternative approach that bypasses line search interventions while assuring monotonicity of the objective as a function of  iteration count. The chopped random basis (CRAB) \cite{doria2011optimal, caneva2011chopped} and gradient optimization of analytic controls (GOAT) \cite{machnes2018tunable} algorithms are more recent alternatives that enable the user to regulate the frequency spectrum of the control signal. These latter methods expand the signal over finitely many Fourier modes and treat the coefficients and frequencies as the optimization variables.  The benefits and limitations of each control algorithm have been compared in terms of their performances on certain quantum systems \cite{machnes2011comparing, riaz2019optimal}.

The parameters of the system are often susceptible to experimental uncertainty in which the assumed values of the control model may vary from those of the plant over a range that generally depends on the application and calibration capability. Furthermore, spatial and temporal variations in the qubit architecture, environment, and control hardware often limit the usability of sensing or computing hardware. Because feedback measurements are generally not plausible in the coherent control of quantum systems, it is often desirable to design open-loop controllers that are robust to parameter uncertainty \cite{koch2022quantum}.  A common approach is based on sampling the uncertainty set at finitely many values and collecting the resulting equations into one deterministic system.  Stochastic gradient-descent \cite{ge2021risk} and supervised machine learning \cite{wu2019learning} over batches of randomly chosen parameter values have also been proposed to further improve robustness over fixed sampling. The GRAPE, Krotov method, CRAB, and GOAT algorithms can be modified in theory to accommodate these robust extensions. Another approach that has been developed for general linear and bilinear systems is the method of moments \cite{zeng2016moment,narayanan2020moment,ning2022nmr,baker2024convergence}.  This method is based on polynomial expansion of the state vector over the uncertainty set in which the expansion coefficients are time-varying vector-valued functions.  In contrast to sampling, the method of moments with orthogonal polynomials offers exponential rates of decay in the norm of the coefficients under certain regulatory and controllability assumptions \cite{wang2012convergence}.

In this study, we develop a robust optimal control algorithm for quantum gate preparation in the presence of decoherence with the environment. The key contribution of this paper is the integration of robustness to parameter inhomogeneity and uncertainty into the evolution equations by introducing external parameters to the Hamiltonian and Lindbladian operators. This enables polynomial expansion of the state over the uncertainty parameter set so that the control design can directly account for variation from nominal model parameters through a small number of polynomials.  The control synthesis is based on adaptive linearization of the quantum evolution operator in addition to an iterative quadratic program (QP) with which the amplitudes of an arbitrarily chosen control signal are gradually shaped between iterations until they converge to an optimal control signal. Quadratic programming provides the user with flexibility to regulate the amplitude and rate of change of the control signal, which may be required in practice.  We prove that if the system is deterministic and the signal restrictions are relaxed, then the proposed control algorithm can be reduced to conventional GRAPE. The algorithm is demonstrated for preparation of quantum Controlled Not (CNOT) and SWAP gates in which the signal is restricted in amplitude and is robust to its intensity and the interaction strength of the qubits.

The rest of the paper is organized as follows.  Section \ref{sec:motion_eqs} defines the Lindbladian evolution of a system of two interacting qubits.  Also in this section, the density matrix and evolution equations are expressed in matrix-vector form and the robust expansion of the vectorized density matrix over a finite number of Legendre polynomials is provided in full detail. The iterative control algorithm is presented in Section \ref{sec:control_algorithm} and a proof of reduction to GRAPE under the above stated relaxations is presented.  Section \ref{sec:examples} showcases the performance of the control algorithm for preparation of two common quantum gates. Conclusions and potential advancements of the algorithm are discussed in Section \ref{sec:conc}.

\section{Evolution of Open Qubit Systems} \label{sec:motion_eqs}

An open qubit system evolving under the influence of electromagnetic fields applied locally to each qubit is represented in terms of the dynamics of its density matrix.  The dynamical equations are vectorized and parameters are introduced into the resulting equations to compensate for uncertainty in nominal values internal to the system.  This results in an uncountable collection of dynamical systems indexed by the introduced uncertain parameters.  The state of the collection of dynamical systems is then expanded over finitely many Legendre polynomials to produce a finite-dimensional approximation.

\subsection{The Lindblad Equation} \label{sec:lindblad}

For an isolated system, the quantum dynamical evolution is governed by the Schr{\"o}dinger equation. For a quantum system in contact with an environment, the Lindblad equation provides a way to account for effects driven by energy exchange and the decay of coherence.  The general form of the Lindblad master equation is given by
\begin{equation} \label{eq:Lindblad_motion}
     \frac{d \rho}{dt} = -i [  H, \rho] + \sum_{m
    }  \Gamma_m \rho   \Gamma_m^\dagger - \frac{1}{2} \{   \Gamma_m^\dagger   \Gamma_m, \rho \}  ,
\end{equation}
where $H(t)$ and  $\rho(t)$ are the Hamiltonian and density matrices, respectively, each of size $d\times d$, and $i$ is the imaginary unit.
The unit of energy (Hamiltonian) and time $t$ are chosen so that $\hbar = 1$. 
The Lindblad matrices $\Gamma_k$ encode the effect of the environment on the quantum system. These matrices will be defined below. 
We follow convention and denote by $[A,B]=AB-BA$ and $\{A,B\}=AB+BA$ the commutation and anti-commutation operations of square matrices $A$ and $B$. 
Moreover, the complex conjugate transpose of a matrix $A$ is denoted by $A^{\dagger}$. 
The commutation term in the Lindblad equation describes the unitary evolution governed by the Hamiltonian. The other terms integrate non-unitary effects and energy exchange with the environment by way of dissipation and decoherence.
The particular form of the Lindblad equation ensures trace preserving evolution that results in a physically relevant Hermitian density matrix.

Consider an interacting system of two qubits exposed to time-varying electromagentic fields. In a rotating reference frame, the Hamiltonian can be written \cite{lu2017enhancing,vstelmachovivc2004quantum,nigmatullin2009implementation}
\begin{equation} \label{eq:Hamiltonian}
        H=\frac{J}{4}\sigma_1^z\sigma_2^z+\frac{1}{2}\left(u_1^x\sigma_1^x+u_1^y\sigma_1^y+u_2^x\sigma_2^x+u_2^y\sigma_2^y\right),
\end{equation}
where the matrices $\sigma^{m}_{1}=\sigma^m\otimes I_2$ generate rotations on the first qubit about the $m=x$, $y$, and $z$ axes, respectively, and likewise for $\sigma^m_{2}=I_2\otimes \sigma^m$ on the second qubit. Here, $I_2$ is the $2\times 2$ identity matrix, $\sigma^m$ are the Pauli matrices for $m=x$, $y$, and $z$, and $A\otimes B$ is the Kronecker product of matrices $A$ and $B$. The above Hamiltonian represents an Ising interaction of strength $J$ in which electromagnetic fields of amplitudes $u_n^x(t)$ and $u_n^y(t)$ are locally applied on the $n$-th qubit. 
The Linblad matrices appearing in Eq. \eqref{eq:Lindblad_motion} are defined by $ \Gamma_z= \sqrt{\gamma}(\sigma_1^z+\sigma_2^z)/2$ and $  \Gamma_{\pm}= \sqrt{\gamma} (\sigma_1^x+\sigma_2^x)/2\pm i   \sqrt{\gamma}(\sigma_1^y+\sigma_2^y)/2$, where $\gamma$ is a uniform rate of decoherence. It is straightforward to verify using the Lie rank test \cite{d2021introduction} that the above Ising interaction of two qubits with two local controls for each qubit is completely controllable in the case of closed systems represented by zero decoherence rates. If decoherence rates are nonzero, then the evolution equations are generally not fully controllable \cite{koch2016controlling}. We examine the degradation of controllability below with numerical examples.

\subsection{Vectorizing the Density Matrix} \label{sec:vectorizing}

To obtain a vector representation of the $d\times d$ density matrix $\rho$, we define the column vector $\vec R=[\rho_{1,1},\rho_{2,1},\dots,\rho_{d,1},\rho_{1,2},\dots,\rho_{d,d}]'$, where $A'$ denotes the transpose of a vector or matrix $A$. The components of the vector $\vec R$ are decomposed into real and imaginary parts to form the real column vector $\vec{\mathcal X} = [\text{Re}(\vec R'),\text{Im}(\vec R')]'$, where real and imaginary operations on the vector $\vec R$ are applied componentwise.   Define the real square matrices
\begin{eqnarray*}
\tilde{\mathcal B}_n^x &=& \frac{1}{2}\left(I_d\otimes   \sigma_n^x -  \left(\sigma_n^x\right)'\otimes I_d\right), \\
 \tilde{\mathcal B}_n^y &=& -\frac{i}{2}\left(I_d\otimes  \sigma_n^y -  \left(\sigma_n^y\right)'\otimes I_d\right), \\
   \tilde{\mathcal A}^{L} &=& \sum_{m}   \Gamma_m\otimes   \Gamma_m-\frac{1}{2}\left(I_d\otimes   \Gamma_m'  \Gamma_m+   \Gamma_m'  \Gamma_m\otimes I_d\right), \\
   \tilde{\mathcal A}^z &=& \frac{J}{4}\left(I_d\otimes  \sigma_1^z\sigma_2^z-  \left(\sigma_1^z\sigma_2^z\right)'\otimes I_d\right),
\end{eqnarray*} 
where $I_d$ is the $d\times d$ identity matrix. With these definitions, the evolution of the vector representation $\vec{\mathcal X} $ of $\rho$ is given by
\begin{equation} \label{eq:dynamics_vectorized}
    \frac{d \vec{\mathcal X}}{dt}  =\alpha \mathcal A\vec{\mathcal X} +\beta \left(u_1^x\mathcal B_1^x +u_1^y \mathcal B_1^y+u_2^x\mathcal B_2^x +u_2^y \mathcal B_2^y \right) \vec{\mathcal X},
\end{equation}
with block matrices defined by $\mathcal B_n^y=\text{diag}(\tilde{\mathcal B}_n^y,\tilde{\mathcal B}_n^y)$ and 
\begin{equation}
    \mathcal A=\begin{bmatrix}
        \tilde{\mathcal A}^{L} & - \tilde{\mathcal A}^z \\
        \tilde{\mathcal A}^z &\tilde{\mathcal A}^{L}
    \end{bmatrix},
    \quad 
    \mathcal B_n^x=\begin{bmatrix}
        0 & -\tilde{\mathcal B}_n^x \\
        \tilde{\mathcal B}_n^x & 0
    \end{bmatrix}.
\end{equation}
We have introduced parameters $\alpha\in [\alpha_{\min}, \alpha_{\max}]$ and $\beta\in [\beta_{\min}, \beta_{\max}]$ to compensate for uncertainty in the assumed nominal values of the internal parameters.  This results in a state vector $\vec{\mathcal X}(t;\alpha,\beta)$ governed by an uncountable collection of structurally identical dynamical systems parameterized by $\alpha$ and $\beta$.  Our goal is to synthesize a single open loop control signal for $u_n^x(t)$ and $u_n^y(t)$ that actuate a robust quantum gate transition independently of the values of $\alpha$ and $\beta$ in their respective compact intervals of uncertainty. 


\subsection{Integrating the Robust Parameters}

The intervals $[\alpha_{\min},\alpha_{\max}]$ and $[\beta_{\min},\beta_{\max}]$ over which the uncertain parameters $\alpha$ and $\beta$ are defined will be shifted to the domain $[-1,1]$ of the Legendre polynomials.  The transformations are defined as
\begin{equation} \label{eq:domain_transformations}
   \alpha(a)=\underline \alpha a+\overline \alpha, 
   \qquad  \beta(b)=\underline \beta b+\overline \beta, 
\end{equation}
where we use the notation $\overline \xi=(\xi_{\max}+ \xi_{\min})/2$ and $\underline \xi=(\xi_{\max}- \xi_{\min})/2$ for $\xi=\alpha$ and $\xi=\beta$.  Using completeness of the Legendre polynomials, we write
\begin{equation} \label{eq:Leg_expansion}
    \vec X(t;\alpha(a),\beta(b))=\sum_{p,q=0}^{\infty} \vec x_{p,q}(t)L_p(a)L_q(b),
\end{equation}
where $L_q$ is the normalized Legendre polynomial and the expansion coefficients are defined by
\begin{equation} \label{eq:expansion_coefficients}
    \vec x_{p,q}(t)=\int_{-1}^1 \int_{-1}^{1}\vec X(t;\alpha(a),\beta(b)) L_p(a)L_q(b) da db.
\end{equation}
Note that each coefficient $\vec x_{p,q}$ is a vector of the same size as $\vec X$. By truncating the expansion, we obtain
\begin{equation} \label{eq:truncated_series}
    \vec X(t;\alpha(a),\beta(b))\approx\sum_{p=0}^{N_{\alpha}}\sum_{q=0}^{N_{\beta}}\vec x_{p,q}(t)L_p(a)L_q(b),
\end{equation}
where $N_{\alpha}$ and $N_{\beta}$ are the maximum degrees of the Legendre polynomials defined on the respective parameter intervals. These degrees are generally dependent on the parameters of the problem and the desired target fidelity.

By differentiating \eqref{eq:expansion_coefficients}, substituting the time-derivative of $\vec X$ according to \eqref{eq:dynamics_vectorized} into the resulting expression, and applying the recurrence relation of Legendre polynomials defined by $
    \xi L_p(\xi)=c_{p-1}L_{p-1}(\xi)+c_pL_{p+1}(\xi)$, where $c_p=(p+1)/{\sqrt{(2p+3)(2p+1)}}$, we obtain the necessary dynamics of the coefficients.  Define the column vector of stacked coefficients of the form $\vec x=[\vec x_{0,0}',\dots, \vec x_{0,N_{\beta}}',\dots,\vec x_{N_{\alpha},0}',\dots,\vec x_{N_{\alpha},N_{\beta}}']'$.  The desired dynamics of $\vec x$ are finally obtained as
    \begin{equation} \label{eq:bilinear_ODE}
        \frac{d\vec x}{dt}= A \vec x +\left(u_1^x B_1^x +u_1^y B_1^y+u_2^x B_2^x +u_2^yB_2^y \right) \vec x,
\end{equation}
where the state and control matrices are defined by $A= C_{\alpha} \otimes I_{N_{\beta}+1} \otimes \mathcal A $ and $B^{m}_{n}=  I_{N_{\alpha}+1} \otimes C_{\beta} \otimes \mathcal B^m_{n}$ for $m=x,\;y$ and $n=1,\;2$. 
Here, the tri-diagonal symmetric matrices $C_{\xi}$ are defined by
\begin{IEEEeqnarray*}{lll}
    C_{\xi}=\begin{bmatrix}
        \overline \xi & c_0\underline \xi &  &  & &    \\
        c_0\underline \xi & \overline \xi & c_1\underline \xi &  & &  \\
         & c_1\underline \xi & \overline \xi &  & &   \\
         & & & \ddots & &  \\
         & & & & \overline \xi &  c_{N_{\xi}-1}\underline \xi \\
         & & & & c_{N_{\xi}-1}\underline \xi &  \overline \xi
    \end{bmatrix}
\end{IEEEeqnarray*}
for $\xi=\alpha$ and $\xi=\beta$. By solving for $\vec x$, we obtain the solution $\vec X$ using Eq. \eqref{eq:Leg_expansion}. Our subsequent exposition is done for the finite-dimensional system in Eq. \eqref{eq:bilinear_ODE}. We conclude this section by referencing a potential avenue for alternative model formulations. The products of Legendre polynomials appearing in Eq. \eqref{eq:Leg_expansion} could be replaced with products of the form $P_p(\alpha)Q_q(\beta)$, where $P$-s and $Q$-s represent functions from distinct orthonormal basis sets.  For example, a specific problem may require periodicity over $\beta$, from which the user may wish to expand the state over products of Legendre polynomials in $\alpha$ and sinusoids in $\beta$.  The derivations above could be tailored to accommodate such alternative formulations.  

\section{Robust Open-Loop Control} \label{sec:control_algorithm}

We now formulate a procedure to synthesize open-loop controls for robust quantum gate preparation and propose a numerical algorithm to obtain a solution.  The algorithm is based on iterative quadratic programming, by which an arbitrary control signal is iteratively perturbed until converging to an optimal signal. Each perturbation is approximated with a time-varying linear model that is updated after each iteration to adequately represent the local nonlinear dynamics. 

\subsection{Control Formulation}
The qubit system in Eq. \eqref{eq:bilinear_ODE} is initially in a state denoted by $\vec x_0$, so that
\begin{eqnarray} \label{eq:initial_cond}
    \vec x(0) =\vec x_0.
\end{eqnarray}
Given a desired target state $\vec x_T$, we seek open-loop control signals 
 that maximize the standard inner product $\langle \vec x_T|\vec x(T)\rangle$ or minimize the error induced in the Euclidean norm, where $T$ represents the elapsed time at which the transition to the target state should be achieved.   Because $\vec x_T$ and $\vec x(T)$ are vectors with real components, formulations that optimize the inner product and error formulation objectives are equivalent in theory.  We implement both approaches for completeness, flexibility, and to analyze the effects these formulations have on the iterative algorithm.

To ensure a physically relevant outcome, restrictions on the amplitude of the signal and its time-derivative are enforced with inequality constraints of the form
\begin{equation} \label{eq:control_inequalities}
    u_{\min}\le u_{n}^{m}(t)\le u_{\max}, \; \Delta u_{\min}\le \dot u_{n}^{m}(t) \le \Delta u_{\max}.
\end{equation}
These must be satisfied simultaneously for both $m=x,\;y$, both $n=1,\;2$, and  all $t\in [0,T]$.  The amplitude bounds $u_{\min}$ and $u_{\max}$ and derivative limits $\Delta u_{\min}$ and $\Delta u_{\max}$ are parameters specified by the user. The optimal control problems (OCPs) for maximizing the inner product or minimizing the error are formulated as
\begin{equation} \label{eq:ocp}
\begin{array}{ll}
\text{minimize:}   & -\langle \vec x_T|\vec x(T)\rangle \quad \text{or} \quad \|\vec x(T)-\vec x_T\|^2 \\
\text{subject to:}   &\text{Robust evolution in Eq. \eqref{eq:bilinear_ODE}}, \\
&\text{Initial condition in Eq. \eqref{eq:initial_cond}}, \\
 &\text{Signal restrictions in Eq.  \eqref{eq:control_inequalities}}.
 \end{array}
\end{equation}
Because these problems contain quadratic equality constraints, analytical solutions are only attainable for select examples and numerical algorithms are generally necessary.

\subsection{Time Discretization and Linearization}

The time interval $[0,T]$ is discretized into $K$ equal steps $[t_{k},t_{k+1})$ of duration $\Delta t=T/K$. In the event that the control signal is constant during each step, the bilinear system in Eq. \eqref{eq:bilinear_ODE} reduces to a linear system by which the state evolves according to the relation $\vec x_{k+1}=U_k\vec x_k$,
where $U_k$ is the matrix exponential defined by
\begin{equation} \label{eq:unitary_matrix}
    U_k=\text{exp} \; \Delta t \left(A+\sum_{n=1}^2 u^x_n(t_k)B^x_n+u^y_n(t_k)B^y_n\right).
\end{equation}
Subscript notation $\vec x_k$ is used to denote the vector sample $\vec x(t_k)$ and similarly for matrices.  By repeatedly applying the transition one time step in succession, we obtain a representation between initial and terminal states given by
\begin{equation}  \label{eq:unitary_evolution}
    \vec x_K=U_{K-1}U_{K-2}\cdots U_1U_0 \vec x_0.
\end{equation}
Because the individual time samples of the control signals are the decision variables of the OCP together with the fact that they are contained in the matrix exponentials, the above transition between the initial and terminal states is highly nonlinear as a function of the control signal.  In order to use efficient quadratic programs, we consider local linear dynamics of the terminal state $\vec x_K$ that result from slight perturbations of a defined control signal.

Differentiating $\vec x_K$ in Eq. \eqref{eq:unitary_evolution} with respect to the individual time sample $u^m_n(t_k)$ of the control signal $u_n^m$ gives the expression
\begin{equation} \label{eq:unitary_evolution_derivative}
   \frac{\partial \vec x_K}{\partial u^m_n(t_k)}=U_{K-1}\cdots U_{k+1}\frac{\partial U_k}{\partial u^m_n(t_k)}U_{k-1}\cdots U_0 \vec x_0.
\end{equation}
We follow the conventional GRAPE algorithm \cite{khaneja2005optimal} and apply the second order approximation
\begin{eqnarray} \label{eq:approx_derivative}
    \frac{\partial U_k}{\partial u^m_n(t_k)} \approx \Delta t B^m_nU_{k}
\end{eqnarray}
to each of the four control signals $u_n^m$ ($m=x,y$ and $n=1,2$) at every time index $k=0,1,\dots,K-1$.  For each $k$, define $\vec u_k=[u_1^x(t_k),u_1^y(t_k),u_2^x(t_k),u_2^y(t_k)]'$ and use these vector time samples to construct the $4K$-dimensional control vector $\vec u=[\vec u_0',\vec u_1',\dots,\vec u_{K-1}']'$.  Likewise, define the matrices
\begin{equation}
    Q_k=\left[ \frac{\partial \vec x_K}{\partial u^x_1(t_k)},\frac{\partial \vec x_K}{\partial u^y_1(t_k)},\frac{\partial \vec x_K}{\partial u^x_2(t_k)}, \frac{\partial \vec x_K}{\partial u^y_2(t_k)} \right]
\end{equation}
and construct the Jacobian $Q=[Q_0,Q_1,\dots,Q_{K-1}]$. Given $\vec u$ and the corresponding terminal state $\vec x_K$ computed using Eq. \eqref{eq:unitary_evolution}, we wish to find an update $\vec u \rightarrow \vec u+\delta \vec u$ to the signal and a corresponding update $\vec x_K\rightarrow \vec x_K+\delta \vec x_K$ to the terminal state that improves the inner product or induced error between the terminal and target states. To first order, the perturbations are related through the linear representation
\begin{eqnarray} \label{eq:linear_transition}
    \delta \vec x_K=Q\delta \vec u.
\end{eqnarray}
We are now in position to define the control algorithm.

\subsection{Iterative Quadratic Programming}

A perturbation $\delta \vec u$ in the control signal $\vec u$ that seeks to improve the objective in Eq. \eqref{eq:ocp} is determined by maximizing the inner product according to the QP given by
\begin{equation} \label{eq:quadratic_program_inner_prod}
\begin{array}{ll}
\text{maximize:}   & \langle \vec x_T|Q \delta \vec u\rangle-\frac{1}{2}\lambda \|\delta \vec u \|^2 \\
\text{subject to:}   &\text{Signal restrictions in Eq.  \eqref{eq:control_inequalities}},
 \end{array}
\end{equation}
or by minimizing the induced error according to
\begin{equation} \label{eq:quadratic_program_error}
\begin{array}{ll}
\text{minimize:}   & \frac{1}{2}\|Q\delta \vec u +\vec x_K-\vec x_T\|^2+\frac{1}{2}\lambda \|\delta \vec u \|^2 \\
\text{subject to:}   &\text{Signal restrictions in Eq.  \eqref{eq:control_inequalities}}.
 \end{array}
\end{equation}
A positive regulation parameter $\lambda$ is introduced into each program to restrict the size of the perturbation between iterations so that the local dynamics and hence the accuracy of linearization can be regulated to arbitrary precision.  The selected values for $\lambda$ are generally different for the two programs.  The solution $\delta \vec u$ of either quadratic program is used to update the control $\vec u \rightarrow \vec u+\delta \vec u$ with which the corresponding terminal state $\vec x_K$ and transition matrix $Q$ are subsequently updated accordingly.  These updates are passed back to the respective optimization problem, the parameter $\lambda$ is subsequently decreased if the updated objective in Eq. \eqref{eq:ocp} decreased, and the algorithm is repeated until a convergence tolerance is satisfied.

We conclude this section by analyzing the iterative update when signal restrictions are relaxed.  In the absence of signal restrictions, the optimization problems \eqref{eq:quadratic_program_inner_prod} and \eqref{eq:quadratic_program_error} each reduce to unconstrained QPs for which exact solutions can be derived.  In the case of maximizing the inner product, the solution is given by
\begin{eqnarray}
    \delta \vec u = \frac{1}{\lambda}Q'\vec x_T.
\end{eqnarray}
This expression is interpreted as gradient ascent with approximate derivatives \eqref{eq:approx_derivative} and a line search parameter $1/\lambda$.  If robust control is bypassed by writing $N_{\alpha}=0$ and $N_{\beta}=0$ in Eq. \eqref{eq:truncated_series}, then the above expression with appropriate step size is equivalent to gradient ascent obtained by conventional GRAPE.  Let us turn our attention to minimizing the induced error. If signal restrictions are relaxed, then the solution is given by
\begin{eqnarray}
    \delta \vec u = -\left(Q'Q+\lambda I\right)^{-1}Q'\left(\vec x_K-\vec x_T \right),
\end{eqnarray}
where $I$ is the $4K \times 4K$ identity matrix.  This is the famous Levenberg update for the iterative solution of nonlinear equations \cite{levenberg1944method}.  Unlike gradient ascent, the regulation parameter in the Levenberg update influences both the length and direction of the update. In either case, it can be verified that $\delta \vec u \rightarrow \vec 0$ as $\lambda \rightarrow \infty$. This confirms that the accuracy of the linear system approximation may be regulated to arbitrary precision between iterations by choosing $\lambda$ sufficiently large. We use the formula $\lambda=\lambda_0\|\vec x_K-\vec x_T\|^2$ to achieve quadratic convergence \cite{yamashita2001rate}, where $\lambda_0$ is a constant that doubles adaptively if the update worsens the error from the previous iterate.  To summarize, we have the following result.
\begin{prop}
Suppose that (i) the robustness portion of the algorithm is bypassed by choosing $N_{\alpha}=N_{\beta}=0$ in Eq. \eqref{eq:truncated_series} and (ii) signal restrictions are removed from the quadratic programs \eqref{eq:quadratic_program_inner_prod}-\eqref{eq:quadratic_program_error}.  Then the solutions provided by \eqref{eq:quadratic_program_inner_prod} and \eqref{eq:quadratic_program_error} are equal to gradient and Levenberg updates, respectively.
\end{prop}

\section{Quantum Gate Preparation Examples}
\label{sec:examples}
We now demonstrate the robust optimal control algorithm for two important two-qubit gates (CNOT \cite{zajac2018resonantly} and SWAP \cite{ono2017implementation}). 
We use conventional notation to denote the state of a single qubit and that of a system of two interacting qubits.  For a single qubit, an arbitrary state $|A\rangle$ is identified as a two-dimensional vector with complex components spanned by the basis states $|0\rangle=[1,0]'$ and $|1\rangle=[0,1]'$.  For a two-qubit system, a general state is spanned by product states of the form $|AB\rangle=|A\rangle\otimes |B\rangle$, where $|A\rangle$ and $|B\rangle$ are states of the first and second qubits, respectively.  For the Hamiltonian in Eq. \eqref{eq:Hamiltonian}, let us note that the two-qubit basis states $|00\rangle$, $|01\rangle$, $|10\rangle$, and $|11\rangle$ respectively correspond to the first, second, third, and fourth diagonal components of the density matrix. In the following, all computations are performed in Matlab R2023a on a MacBook Pro with 32 GB of usable memory and an Apple M2 Max processing chip. The general-purpose Matlab function \verb+quadprog+ is used at each iteration with the sparse-linear-convex algorithm and sparse linear algebra. 

\subsection{CNOT Gate}

The CNOT gate is an important component of quantum computation because of its ability to generate entanglement between two qubits. Indeed, the CNOT gate in combination with single qubit gates collectively form a universal set of gates for quantum computation. The action of the CNOT gate 
is defined by the interchange of states $|10\rangle \leftrightarrow |11\rangle$ while not affecting states $|00\rangle$ and $|01\rangle$.  To demonstrate the capabilities of the control algorithm, we design a robust CNOT gate that compensates for 100\% uncertainty in the qubit interaction strength $J$ by choosing the domain $[0,2]$ for the uncertainty parameter $\alpha$. Furthermore, the design simultaneously compensates for 20\% uncertainty in the intensity of the exerted signal by choosing the interval $[0.8,1.2]$ for $\beta$.  We use the parameters $J=0.1$, $T=1$, $K=100$, $N_{\alpha}=1$, and $N_{\beta}=2$, and restrict the amplitudes of the control signal to the interval $[-10,10]$.  Because the complex vector representation $\vec R$ of $\rho$ has 16 components, the size of the Jacobian $Q$ for the robust system in Eq. \eqref{eq:bilinear_ODE} is $32(N_{\alpha}+1)(N_{\beta}+1)\times mK=192\times 400$.

\begin{figure*}[t]
    \centering
    \begin{subfigure}[t]{0.5\textwidth}
        \centering
        \includegraphics[width=1\linewidth]{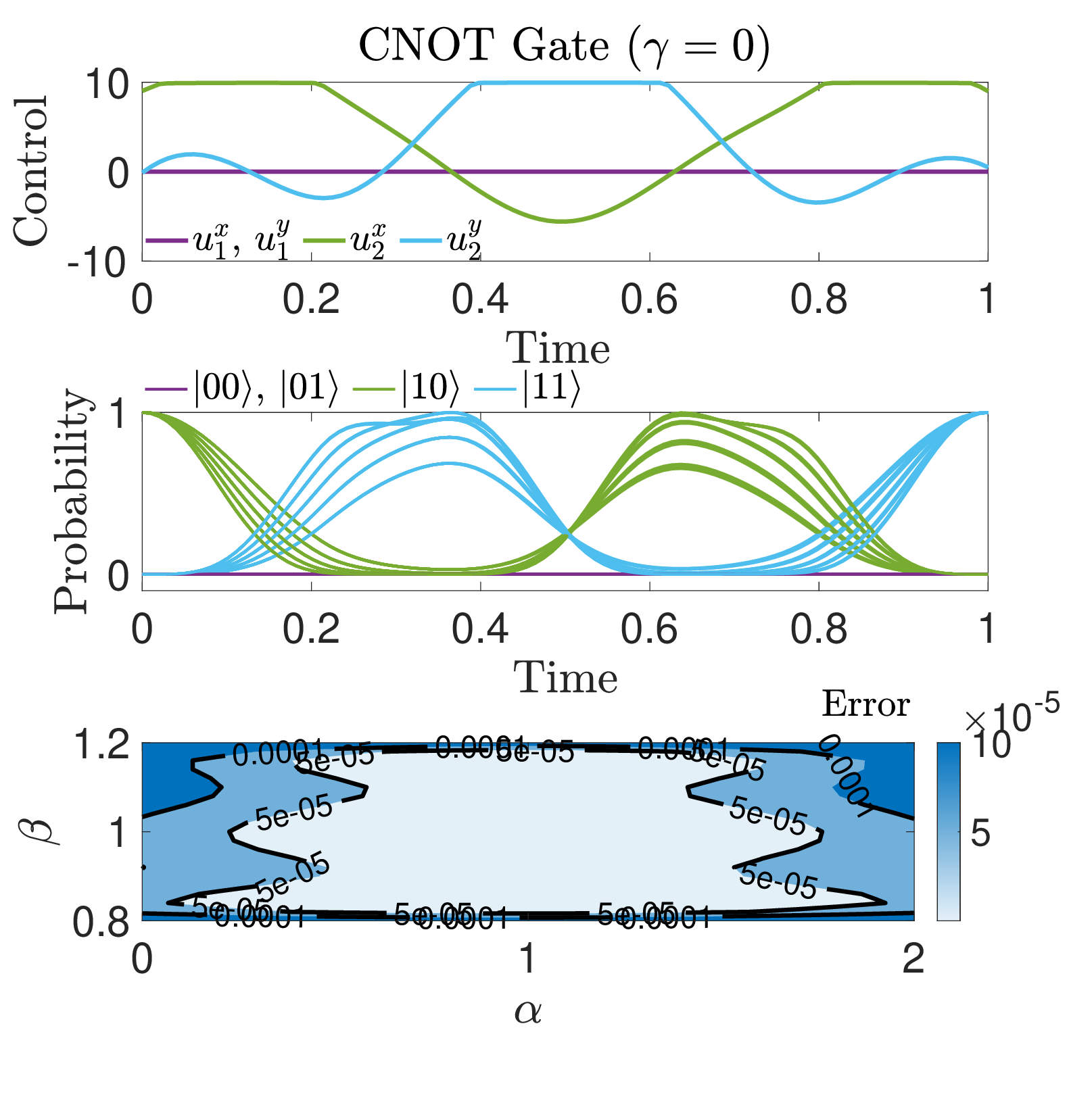}
    \end{subfigure}%
    ~ 
    \begin{subfigure}[t]{0.5\textwidth}
        \centering
        \includegraphics[width=1\linewidth]{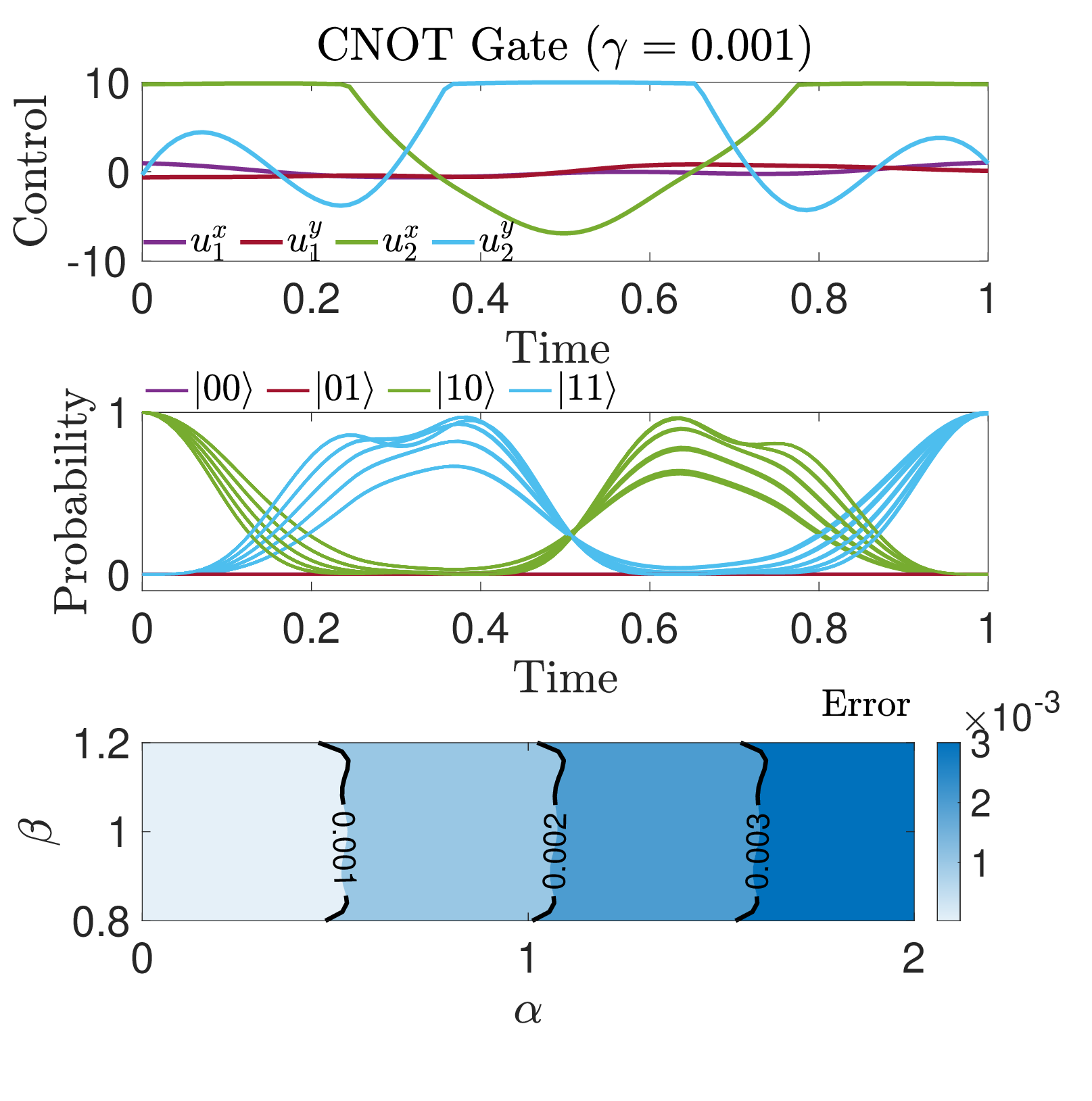}
    \end{subfigure}
    \vspace{-6ex}
\caption{Robust CNOT gate preparation for the two-qubit system.  Depicted from top to bottom are control signals, probabilities of states, and error contours over the region of uncertainty for (left column) zero interaction with the environment and (right column) interation of strength $\gamma=0.001$.}
\label{fig:cnot}
\end{figure*}

\begin{figure*}[t]
    \centering
    \begin{subfigure}[t]{0.5\textwidth}
        \centering
        \includegraphics[width=1\linewidth]{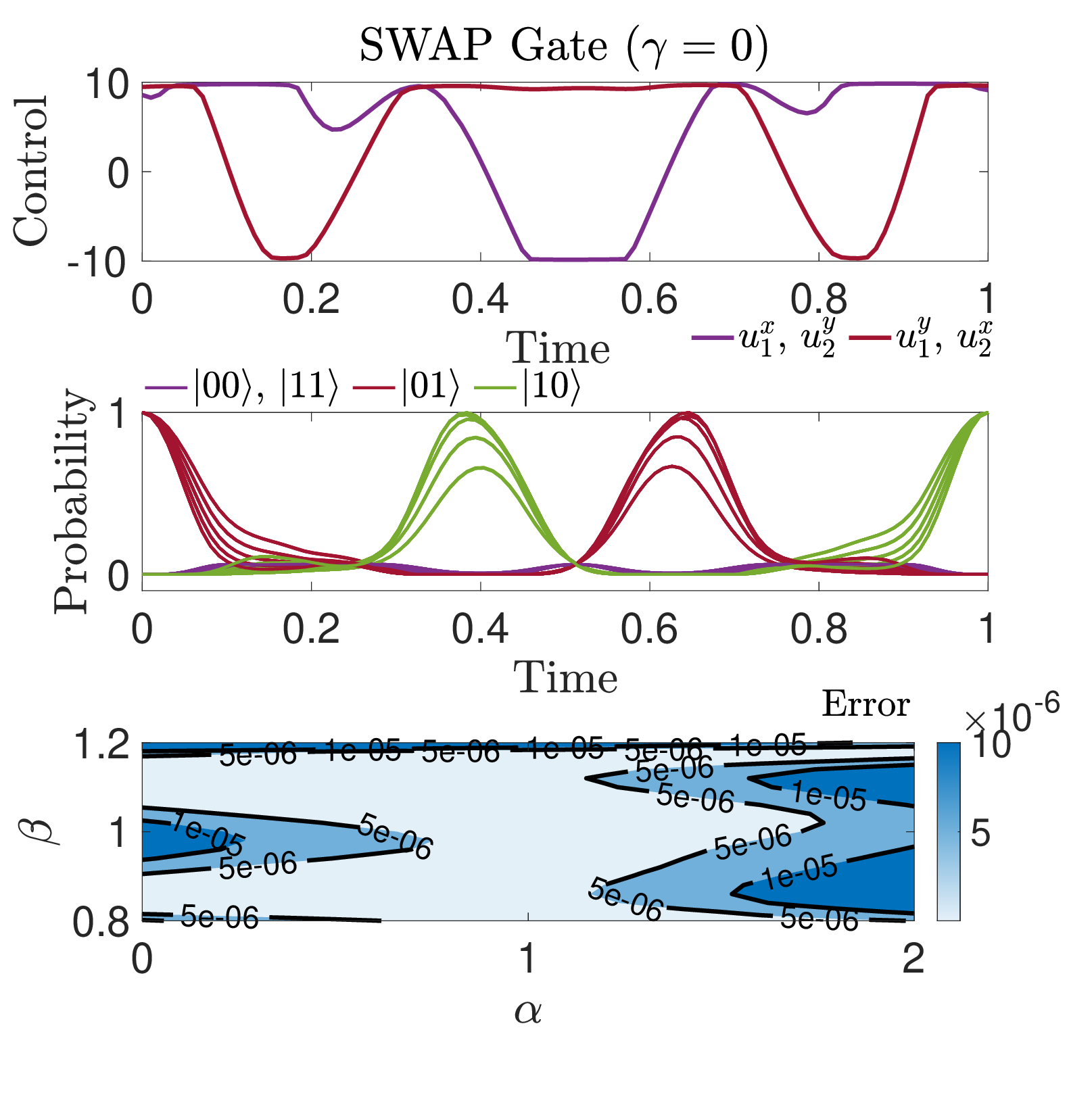}
    \end{subfigure}%
    ~ 
    \begin{subfigure}[t]{0.5\textwidth}
        \centering
        \includegraphics[width=1\linewidth]{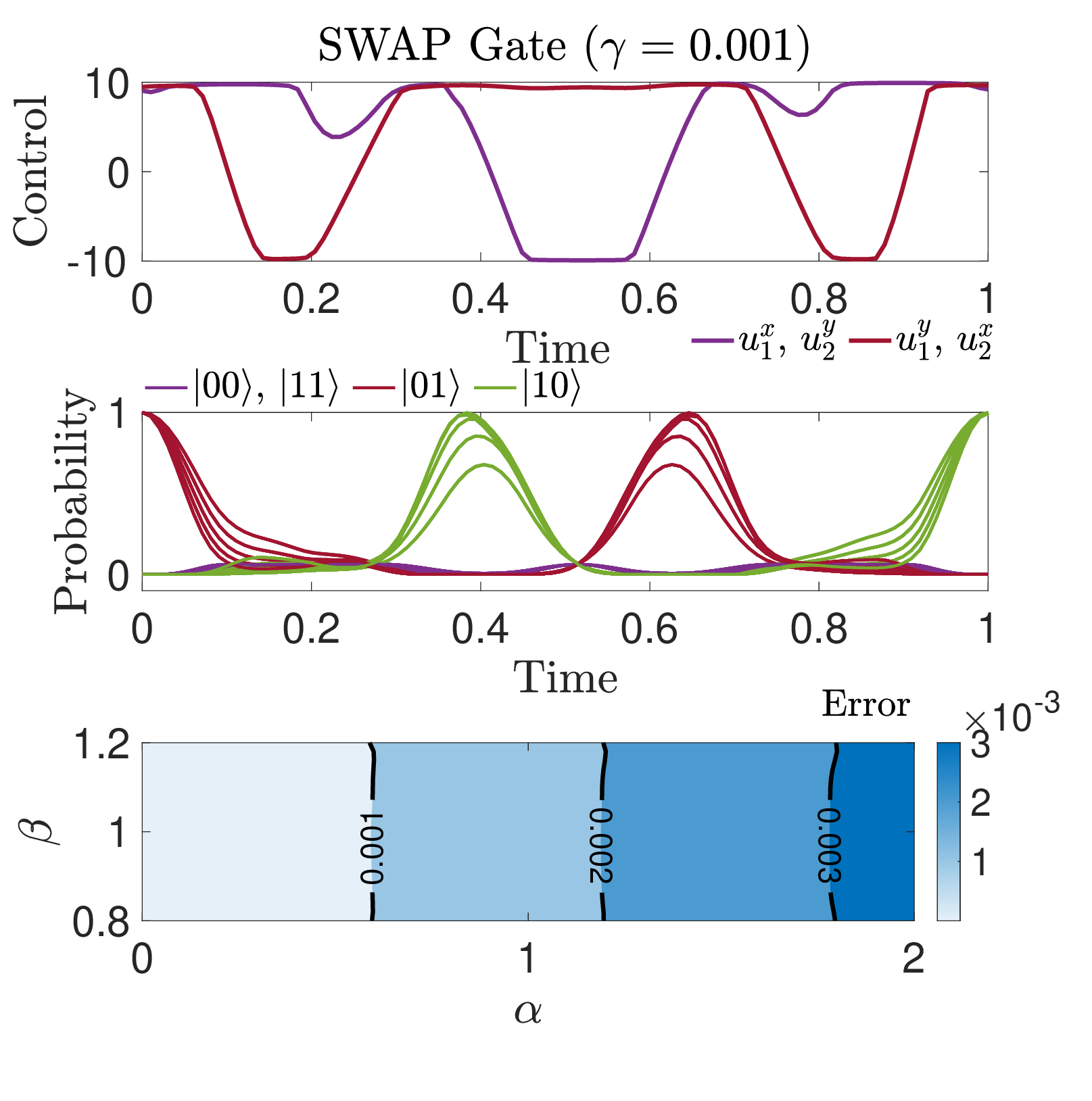}
    \end{subfigure}
        \vspace{-6ex}
\caption{Robust SWAP gate preparation for the two-qubit system.}
\label{fig:swap}
\end{figure*}

The control signals, probability evolution of states, and contours of the terminal error obtained by repeated simulation over the design region of parameters $\alpha$ and $\beta$ for CNOT gate preparation are depicted from top to bottom in Fig. \ref{fig:cnot}. The left column corresponds to zero rates of decoherence and the right column considers the effect of decoherence with a rate equal to $\gamma=0.001$ \cite{riaz2019optimal}. The control algorithm terminates in less than one minute for both problems after reaching the specified number of allowed iterations, which we define to be 50.  It is evident from Fig. \ref{fig:cnot} that the control algorithm successfully prepares a CNOT gate with unprecedented error margins over the region of uncertainty, even in the presence of decoherence.  Although each of the four control signals are present in the control design and are permitted to any amplitudes in the restricted region, the control algorithm for zero decoherence converges to an intuitive signal in which the local controls on the first qubit are zero for all time.  However, as one would expect, in the case with decoherence effects from the environment for which results are shown in the right column, the control algorithm converges to a signal that has nonzero amplitudes applied on the first qubit to mitigate the effects from the environment.  Although still negligible for practical purposes, observe that the error margins increase by nearly two orders of magnitude by increasing the rate of decoherence from $\gamma=0$ to $\gamma=0.001$.  This agrees with the loss of controllability in open environments.

\subsection{SWAP Gate}
The SWAP gate is a fundamental quantum gate that swaps the state of two qubits. It is a key element in quantum communication and entanglement distribution, and can be used to optimize the hardware implementation of a quantum circuit by bringing the states of two qubits close to one another. The SWAP gate acts by interchanging states $|01\rangle \leftrightarrow |10\rangle$, while leaving the states $|00\rangle$ and $|11\rangle$ unaffected.  In the control algorithm, the system parameters stated above are also used for robust SWAP gate preparation.

The control signals, probability evolution of states, and contours of the terminal error over the region of parameters $\alpha$ and $\beta$ for SWAP gate preparation are depicted from top to bottom in Fig. \ref{fig:swap}. As with the CNOT gate, the control algorithm successfully prepares the SWAP gate with similarly minor error margins over the region of uncertainty.  As one may expect, the algorithm determines a symmetric action of the control signals between the two quibits by which we mean that the rotations about the $x$ and $y$ axes generated from $u^x_1$ and $u^y_2$ are equivalent and likewise for $u^y_1$ and $u^x_2$, although all four control signals are synthesized without such symmetry constraints.  In contrast to the states $|00\rangle$ and $|01\rangle$ that remain unaffected during the entire duration of the CNOT gate, the states $|00\rangle$ and $|11\rangle$ that are designed to be unaffected do have positive and nontrivial probability during the evolution until the terminal state when they take on the expected negligible probabilities. However, this is expected because all four control signals are activated during the evolution of the SWAP gate, whereas the control algorithm for CNOT gate preparation determined only two nonzero control signals that were applied locally on the one qubit. 


\section{Conclusion} \label{sec:conc}

This paper develops a robust optimal control algorithm for qubit gate preparation in an open system and successfully demonstrates the capability to prepare highly robust Controlled NOT and SWAP gates. Robustness over a large region of parameter uncertainty with only a small extension of the size of the underlying quantum system is accomplished with the use of a small number of orthogonal polynomials to represent the ensemble state on the domain of the uncertain parameter.  The algorithm utilizes an expansion of the quantum state over the parameterized region of uncertainty with Legendre polynomials that are well-known to provide exponential rates of convergence in the number of polynomials under some regulatory assumptions \cite{wang2012convergence}.  However, we note that other sets of basis functions defined on compact intervals may be considered as a substitute of Legendre polynomials.  In any case, after extending the underlying quantum system to a robust expansion, the control algorithm is then applied to the latter to achieve the desired robust properties. The control algorithm that we have proposed is based on adaptive linearization of the evolution of the qubit system and iterative quadratic programming.  When signal restrictions are removed, the algorithm reduces to conventional GRAPE.

Because of the exponential rate of decay in the coefficients of orthogonal polynomial expansions, the proposed algorithm can be used to improve the robustness of current state-of-the-art algorithms for quantum state preparation, which are based on direct or stochastic sampling \cite{koch2022quantum, ge2021risk, wu2019learning}.  However, there are several avenues to extend and further improve the computation of this polynomial approach.  First, because the polynomial expansion leads to a robust system with the same bilinear structure as the underlying deterministic quantum system, the robust polynomial approach could be integrated with well-developed GRAPE, CRAB, GOAT, and Krylov packages.  In terms of directly improving the proposed algorithm, upper bounds on the difference in the objective function values between iterations could be derived \cite{vu2023iterative} and potentially used for adaptive regulation and monotone convergence. Third, although the robust algorithm uses a small extension of the underlying system, preconditioning and BFGS updates may be used to further improve the speed at which iterations are performed. These extensions could enhance the capabilities of the robust algorithm for implementation on large-scale quantum control problems.

\linespread{1}
\bibliographystyle{unsrt}
\bibliography{references.bib}

\end{document}